\documentclass[12pt,a4paper,reqno]{amsart}
\usepackage{graphics,epic}
\usepackage{amsmath,amssymb, amsthm}
\usepackage[all,2cell]{xy}

\newtheorem{theorem}{Theorem}[section]
\newtheorem*{theorem*}{Theorem}

\newtheorem{lemma}[theorem]{Lemma}
\newtheorem{proposition}[theorem]{Proposition}

\newtheorem*{conjecture*}{Conjecture}

\newtheorem{remark}[theorem]{Remark}

\renewcommand{\hat}[1]{\widehat{#1}}

\newcommand{\m}{\mathfrak{m}}

\newcommand{\ot}{\otimes}

\newcommand{\Z}{\mathbb{Z}}

\newcommand{\C}{\mathbb{C}}

\newcommand{\la}{\lambda}

\def\C{{\mathbb C}}

\def\Z{{\mathbb Z}}

\def\1{{\bf 1}}
\def\la{{\langle}}
\def\ra{{\rangle}}

\def \pf{\noindent {\bf Proof: \,}}
\def\theequation{5.\arabic{equation}}

\def \h{\mathfrak{h}}

\def \w{\omega}
\def \g{\mathfrak{g}}
\def \h{\mathfrak{h}}

\begin{document}

\title[The classification of  extensions of $L_{sl_3}(k,0)$]{The classification of extensions of $L_{sl_3}(k,0)$}
\author{Chunrui Ai}
\address{Chunrui Ai, School of Mathematics and Statistics, Zhengzhou University, Henan 450001, China
}
\email{aicr@zzu.edu.cn}
\author{Xingjun Lin}
\address{Xingjun Lin, Institute of Mathematics, University of Tsukuba, Tsukuba, Japan}
\email{linxingjun88@126.com}
\begin{abstract}
In this paper,  extensions of affine vertex operator algebra $L_{sl_3}(k,0)$, $k\in \Z_+$, are classified by modular invariants.
\end{abstract}
\maketitle
\section{Introduction \label{intro}}
\def\theequation{1.\arabic{equation}}
\setcounter{equation}{0}

One of important problems in the theory of vertex operator algebras is to classify rational vertex operator algebras, and there were many works on classifications. For instance, preunitary vertex operator algebras with central charge $c<1$ were classified in  \cite{DZ}, \cite{M}, \cite{DLin2}; rational  vertex operator algebras with central charge $c=1$ were partially classified in \cite{DM1}, \cite{DJ1}, \cite{DJ2}, \cite{DJ3}; holomorphic vertex operator algebras with central charge 24 were partially classified in \cite{LS}.

It is well-known that vertex operator algebras are closely related to two-dimensional conformal field theory. One of important problems in two-dimensional conformal field theory is to classify modular invariant partition functions, and there were many works on the classification of modular invariant partition functions \cite{CIZ1}, \cite{CIZ2}, \cite{Ka}, \cite{G}. In language of vertex operator algebras, this is equivalent to classify modular invariants of rational vertex operator algebras \cite{DLN}. It was also shown in \cite{DLN} that extensions of rational vertex operator algebras give rise to modular invariants. Therefore, it is interesting to realize modular invariants by constructing extensions of vertex operator algebras explicitly and to classify extensions of rational vertex operator algebras.

In this paper, we consider the classification of extensions of affine vertex operator algebra $L_{sl_3}(k,0)$, $k\in \Z_+$. The method is to classify extensions by modular invariants. This method is first proposed in \cite{KL} in the framework of conformal nets, and was latter used in \cite{DLin2} to classify extensions of $L_{sl_2}(k,0)$, $k\in \Z_+$.  The classification of modular invariants of $L_{sl_3}(k,0)$ has been obtained in  \cite{G}, and plays an important role in this paper. It is shown in this paper that all the extensions of affine vertex operator algebra $L_{sl_3}(k,0)$ can be realized by simple current extensions \cite{DLM4} and conformal embeddings \cite{K}. To complete the classification, the key point is to prove the uniqueness.

The rest of this paper is organized as follows: In Section 2, we recall some basic facts about affine vertex operator algebras. In Section 3, we  classify extensions of $L_{{sl_3}}(k,0)$.

\section{Affine vertex operator algebras $L_{sl_3}(k,0)$}
\def\theequation{3.\arabic{equation}}
\setcounter{equation}{0}
\subsection{Affine vertex operator algebras}\label{subsec1}
 In this subsection, we recall some facts about vertex operator algebras associated to affine Lie algebras from \cite{H}, \cite{FZ}, \cite{LL}. Let $\mathfrak{g}$ be a finite dimensional simple Lie algebra, $\mathfrak{h}$ be a Cartan subalgebra of
$\mathfrak{g}$. We denote the simple roots, simple coroots and fundamental weights of $\mathfrak{g}$ by $\{\alpha_1, \cdots, \alpha_n\}$, $\{\alpha_1^{\vee}, \cdots, \alpha_n^{\vee}\}$ and $\{\Lambda_1,\cdots, \Lambda_n\}$, respectively. We use $h^{(1)}, \cdots, h^{(n)}\in \mathfrak{h}$  to denote the fundamental coweights of $\mathfrak{g}$, i.e.,
$$\alpha_i(h^{(j)})=\delta_{i, j}\ \ \text{for}\ i, j=1, \cdots, n.$$ We shall also use $\theta=\sum_{i=1}^na_i\alpha_i$ to denote the highest root of $\mathfrak{g}$.

 Let $\langle\ ,\ \rangle$ be the normalized Killing form of $\mathfrak{g}$ such that $\langle\theta, \theta \rangle=2$. The corresponding affine Lie algebra is defined as  $\mathfrak{\hat{g}}=\mathfrak{g}\otimes\mathbb{C}[t, t^{-1}]\oplus \mathbb{C}K$  with the commutation relations
$$[x(m), y(n)]=[x, y](m+n)+Km\langle x, y\rangle \delta_{m+n,0},$$
$$[\mathfrak{\hat{g}}, K]=0,$$
for any $x, y\in \g$, $m, n\in \Z$, where $x(m)=x\otimes t^m$.

For a complex number $k$, set
$$V_\mathfrak{g}(k, 0)=\text{Ind}_{\mathfrak{g}\otimes \mathbb{C}[t]\oplus \mathbb{C}K}^{\mathfrak{\hat{g}}}\mathbb{C},$$
where $\C$ is viewed as a $\mathfrak{g}\otimes \mathbb{C}[t]\oplus \mathbb{C}K$-module such that $\mathfrak{g}\otimes \mathbb{C}[t]$ acts as 0 and $K$ acts as $k$. It was  proved in \cite{FZ}  that $V_\mathfrak{g}(k, 0)$ has a vertex operator algebra structure if $k\neq -h^{\vee}$, where $h^{\vee}$ is the dual Coxeter number of $\mathfrak{g}$. Furthermore, it is well-known that the vertex operator algebra $V_\mathfrak{g}(k, 0)$ has a unique maximal proper ideal $\mathcal{J}$ (see \cite{K}). The quotient vertex operator algebra $L_\mathfrak{g}(k, 0)=V_\mathfrak{g}(k, 0)/\mathcal{J}$ is a simple  vertex operator algebra and  is strongly regular if $k\in \mathbb{Z}_+$ (see \cite{FZ}, \cite{DLM1},\cite{DM1}). It is also known that  $$\{L_\mathfrak{g}(k, \lambda)| \lambda\in P^k_{+}\} ,$$
where $P^k_+$ denotes the dominant weight with level $k$ of $\g$ (cf. \cite{K}), is the complete list of nonisomorphic  irreducible $L_\mathfrak{g}(k, 0)$-modules (see \cite{FZ}, \cite{LL}).

\subsection{Modular invariants of affine vertex operator algebra $L_{sl_3}(k,0)$}
We first recall from \cite{DLM3}, \cite{Z} some facts about  modular invariance properties of vertex operator algebras. Let $V$ be a strongly regular vertex operator algebra, $M^0, M^1,\cdots, M^p$ be all the nonisomorphic irreducible $V$-modules. It was proved in \cite{Z} that for any $0\leq i\leq p$, $M^i$  has the form
$$M^i=\bigoplus_{n=0}^{\infty}M^i_{\lambda^i+n}$$
for some $\lambda^i\in \C$ such that $M^i_{\lambda^i}\neq 0$, and $\lambda^i$ is called the {\em conformal weight} of $M^i$. Let $\mathbf{H}=\{\tau\in \mathbb{C}| \text{im}\tau>0\}$.  For any homogeneous element $v\in V$, we define the {\em trace function} associated to $M^i$ as follows:
$$Z_{M^i}(v,\tau): =\text{tr}_{M^i}o(v)q^{L(0)-c/24}=q^{\lambda_i-c/24}\sum_{n\in\Z_{\geq 0}}\text{tr}_{M^i_{\lambda_i+n}}o(v)q^n,$$
where $q=e^{2\pi \sqrt{-1} \tau}$ and $o(v)=v_{\text{wt}v-1}$ is the degree zero operator of $v$ (cf. \cite{Z}). It was proved in \cite{DLM3} and \cite{Z} that $Z_{M^i}(v,\tau)$ converges to a holomorphic function on the domain $|q|<1$.

 Recall that the full modular group $SL(2, \mathbb{Z})$
acts on $\mathbf{H}$ as follows:
$$\gamma:   \tau \mapsto \frac{a\tau+b}{c\tau+d}, ~~ \gamma=\left(\begin{array}{cc}a & b \\c & d \end{array}\right)\in SL(2, \mathbb{Z}).$$
In particular, $SL(2, \mathbb{Z})$ acts on the trace functions. The following result was obtained in \cite{Z} (see also \cite{DLM3}).
\begin{theorem}
Let $V$ be a strongly regular vertex operator algebra, $M^0, \cdots, M^p$ be all the nonisomorphic  irreducible $V$-modules. Then the vector space spanned by $$Z_{M^0}(v, \tau), \cdots, Z_{M^p}(v, \tau)$$ admits a representation $\phi$ of  $SL(2, \mathbb{Z})$ and the transformation matrices are independent of the choice of $v\in V$.
\end{theorem}
\vskip.25cm
Let
$S=\left(\begin{array}{cc}0 & -1  \\1 & 0 \end{array}\right)$ and  $T=\left(\begin{array}{cc}1 & 1 \\0 & 1 \end{array}\right)$. We know that $S, T$ are generators of $SL(2, \mathbb{Z})$. A {\em modular invariant} of $V$ is then defined to be a $(p+1)\times(p+1)$-matrix $X$ satisfying the following conditions:\\
(P1) $X_{0,0}=1$;\\
(P2) The entries of $X$ are nonnegative integers;\\
(P3) $X\phi(S)=\phi(S)X$ and $X\phi(T)=\phi(T)X$.

\vskip.25cm
 It was shown in \cite{DLN} that any {\em extension} $U$ of $V$, i.e., $V$ is a vertex operator subalgebra of $U$ and $V, U$ have the same conformal vector, gives rise to a modular invariant of $V$. More precisely, let $U$ be an extension  of $V$. Assume that $U$ is strongly regular, and let $N^0, \cdots, N^k$ be all the nonisomorphic irreducible $U$-modules. For $u,v\in U$,  set
$$f_U(u,v,\tau_1,\tau_2)=\sum_{i=0}^kZ_{N^i}(u,\tau_1)\overline{Z_{N^i}(v,\tau_2)},$$ where $\tau_1, \tau_2 \in\mathbf{H}$.
Since each irreducible $U$-module $N^i$ is a direct sum of irreducible $V$-modules, there exists a matrix $X=(X_{i,j})$ such that $X_{i,j} \in \Z_{\geq 0}$ for all $i, j$ and  that for any $u,v\in V$,
$$f_U(u,v,\tau_1, \tau_2)=\sum_{i, j=0}^{p}X_{i,j}Z_{M^i}(u, \tau_1)\overline{{Z_{M^j}}(v, \tau_2)}.$$
Moreover, by the following result which was proved in \cite{DLN}, the matrix $X=(X_{i,j})$ is uniquely determined.
\begin{proposition}\label{unique}
Set ${\bf Z}(u,\tau)=(Z_{M^0}(u,\tau), ..., Z_{M^p}(u,\tau))^T.$
 If $A=(a_{ij})$ is a  matrix such that  for any $u,v\in V$, $${\bf Z}(u,\tau_1)^TA\overline{{\bf Z}(v,\tau_2)}=0,$$ then $A=0.$
\end{proposition}

It was further shown in \cite{DLN} that $X$ is a modular invariant of $V$.
\begin{theorem}\label{ccc}
Let $V$, $U$, $X$ be as above. Then $X$ is a modular invariant of $V$.
\end{theorem}
Furthermore, by the discussion above, we have the following simple observation.
\begin{lemma}\label{cre}
Let $V$, $U$, $X$ be as above. Suppose that there exist irreducible $V$-modules $M^i$ and $M^j$ such that $X_{i,j}\neq 0$. Then we have $X_{i,i}\neq 0$ and $X_{j,j}\neq 0$.
\end{lemma}

 \vskip.25cm
 We now recall from \cite{G} some facts about modular invariants of $L_{sl_3}(k,0)$. We use $\rho$ to denote $\Lambda_1+\Lambda_2$,  and $(m, n)$ to denote the weight $\lambda=m\Lambda_1+n\Lambda_2$, where $\Lambda_1, \Lambda_2$ are the fundamental weights of $sl_3$.
Recall that the set of dominant weights with level $k$ is equal to
$$P_+^k=\{m\Lambda_1+n\Lambda_2| m, n\in \mathbb{Z}, 0\leqslant m, n, m+n\leqslant k\}.$$
For convenience, we shall use $P^k$ to denote
$$P_+^k+ \rho=\{m\Lambda_1+n\Lambda_2| m, n\in \mathbb{Z}, 0< m, n, m+n <k+3\}.$$
Define two maps $h: P^k\to P^k $ and $\sigma: P^k\to P^k$ as follows:
\begin{align}\label{120}
h(m, n)=(n, m),\ \ \sigma(m, n)=(k+3-m-n, m).
\end{align}
Hence $\sigma^2(m,n)=(n, k+3-m-n)$.

 Denote the trace function of $L_{sl_3}(k,0)$-module $L_{sl_3}(k,\lambda)$ $ (\lambda=m\Lambda_1+n\Lambda_2\in P_+^k)$ by
$Z_{\lambda+ \rho}^k(u,\tau)=Z_{(m+1,n+1)}^k(u,\tau)$.
Then for any modular invariant $X=(X_{\lambda,\mu})$ ($\lambda, \mu\in P^k$)  of $L_{sl_3}(k,0)$, we will write it in the form
\begin{align*}
\sum_{\lambda, \mu\in P^k}X_{\lambda, \mu}Z_{\lambda}^k(u, \tau_1)\overline{Z_{\mu}^k(v, \tau_2)}.
\end{align*}
It was shown in \cite{G} that
\begin{align*}
\sum_{\lambda, \mu\in P^k}X_{\lambda, h(\mu)}Z_{\lambda}^k(u, \tau_1)\overline{Z_{\mu}^k(v,\tau_2)}
\end{align*}
 also determines a modular invariant $X^C$ of $L_{sl_3}(k,0)$, which is called the {\em conjugate modular invariant} of $X$. In particular, we have $(X^C)_{\lambda, \mu}=X_{\lambda, h(\mu)}$. The following result  was proved in \cite{G} (see also \cite{G1}).
\begin{theorem}\label{aaaa}
 Any modular invariant of $L_{{sl_3}}(k,0)$ is equal to one of the following modular invariants or  their conjugate modular invariants:

{\tiny
$$
\begin{tabular}{|c c c|}
\hline

{\rm any} k  & $\sum\limits_{(m, n)\in P^k}| Z_{(m,n)}^k|^2 $ & $A_{k}$ \\
\hline

$k\not\equiv  0\ ({\rm mod~ 3)}, ~k\geqq 4$ & $\sum\limits_{(m, n)\in P^k}Z_{(m,n)}^k(u, \tau_1)\overline{Z_{\sigma^{k(m-n)}(m,n)}^k(v,\tau_2)}$ & $D_k$ \\
\hline
 $k\equiv 0 \ (\rm mod\ 3)$ & $\frac{1}{3}\sum\limits_{\substack {(m, n)\in P^k \\m\equiv n(\rm {mod}\ 3)}}| Z_{(m, n)}^k + Z_{\sigma(m, n)}^k +Z_{\sigma^2(m, n)}^k|^2$ & $D_k$\\
 \hline
 k=5 &$| Z_{(1,1)}^5+ Z_{(3, 3)}^5|^2 +|Z_{(1, 3)}^5+ Z_{(4, 3)}^5|^2 +|Z_{(3, 1)}^5+ Z_{(3 ,4)}^5|^2$& $ E_5$ \\
 &$+|Z_{(3, 2)}^5+Z_{(1, 6)}^5|^2+ |Z_{(4, 1)}^5+ Z_{(1, 4)}^5|^2+|Z_{(2, 3)}^5+Z_{(6, 1)}^5|^2 $&\\

 \hline
k=9 & $|Z_{(1, 1)}^9+ Z_{(1, 10)}^9+Z_{(10, 1)}^9+Z_{(5, 5)}^9+Z_{(5, 2)}^9+Z_{(2, 5)}^9|^2$ & $ E_9^{(1)}$ \\
&$\ +2|Z_{(3,3)}^9+Z_{(3,6)}^9+Z_{(6,3)}^9|^2$&\\
\hline
k=9 & $|Z_{(1,1)}^9+Z_{(10,1)}^9+Z_{(1,10)}^9|^2+|Z_{(3,3)}^9+Z_{(3,6)}^9+Z_{(6,3)}^9|^2+2|Z_{(4,4)}^9|^2$ & $ E_9^{(2)}$ \\
&$+|Z_{(1,4)}^9+Z_{(7,1)}^9+Z_{(4,7)}^9|^2+|Z_{(4,1)}^9+Z_{(1,7)}^9+Z_{(7,4)}^9|^2+|Z_{(5,5)}^9+Z_{(5,2)}^9+Z_{(2,5)}^9|^2$&\\
&$ \ +\left(Z_{(2,2)}^9(u,\tau_1)+Z_{(2,8)}^9(u,\tau_1)+Z_{(8,2)}^9(u,\tau_1)\right)\overline{Z_{(4,4)}^9(v,\tau_2)} $&\\
&$+Z_{(4,4)}^9(u,\tau_1)\overline{\left(Z_{(2,2)}^9(v,\tau_2)+Z_{(2,8)}^9(v,\tau_2)+Z_{(8,2)}^9(v,\tau_2)\right)}$&\\
\hline
 k=21 & $|Z_{(1,1)}^{21}+Z_{(5,5)}^{21}+Z_{(7,7)}^{21}+Z_{(11,11)}^{21}+Z_{(22,1)}^{21}+Z_{(1,22)}^{21}$ & $ E_{21}$\\
 &$+Z_{(14,5)}^{21}+Z_{(5,14)}^{21}+Z_{(11,2)}^{21}+Z_{(2,11)}^{21}+Z_{(10,7)}^{21}+Z_{(7,10)}^{21}|^2$& \\
 &$ +|Z_{(16,7)}^{21}+Z_{(7,16)}^{21}+Z_{(16,1)}^{21}+Z_{(1,16)}^{21}+Z_{(11,8)}^{21}+Z_{(8,11)}^{21}$&\\
 &$+Z_{(11,5)}^{21}+Z_{(5,11)}^{21}+Z_{(8,5)}^{21}+Z_{(5,8)}^{21}+Z_{(7,1)}^{21}+Z_{(1,7)}^{21}|^2$&\\
 \hline
\end{tabular}
$$}
where we have used $|Z^k_{\lambda_1}+...+Z^k_{\lambda_s}|^2$ $(\lambda_1,...,\lambda_s\in P^k)$ to denote $$(Z^k_{\lambda_1}(u, \tau_1)+...+Z^k_{\lambda_s}(u,\tau_1))\overline{(Z^k_{\lambda_1}(v, \tau_2)+...+Z^k_{\lambda_s}(v,\tau_2))}.$$
\end{theorem}

By a direct calculation, we have the following:
\begin{lemma}\label{modular2}
The modular invariants of types $D_3^C$, $D_6^C$, $(E_9^{(1)})^{C}$, $E_{21}^{C}$  are equal to the modular invariants of types $D_3$, $D_6$, $E_9^{(1)}$,  $E_{21}$, respectively.
\end{lemma}
We also need the following result.
\begin{lemma}\label{5aa}
The modular invariants  $E_5^C$, $(E_9^{(2)})^C$ and   $E_9^{(2)}$ cannot be realized by extensions of $L_{sl_3}(k,0)$.
\end{lemma}
\pf   By a direct calculation, the modular invariants of types $E_5^C$, $(E_9^{(2)})^C$ are equal to
{\tiny
\begin{align*}
&|Z^5_{(1,1)}+Z^5_{(3,3)}|^2+|Z^5_{(4,1)}+Z^5_{(1,4)}|^2+(Z^5_{(1,3)}(u, \tau_1)+ Z^5_{(4,3)}(u, \tau_1))\overline{(Z^5_{(3,1)}(v, \tau_2)+Z^5_{(3,4)}(v, \tau_2))}\\
&+(Z^5_{(3,2)}(u, \tau_1)+Z^5_{(1,6)}(u, \tau_1))\overline{(Z^5_{(2,3)}(v, \tau_2)+Z^5_{(6,1)}(v, \tau_2))}+(Z^5_{(3,1)}(u, \tau_1)+Z^5_{(3,4)}(u, \tau_1))\overline{(Z^5_{(1,3)}(v, \tau_2)+Z^5_{(4,3)}(v, \tau_2))}\\
&+(Z^5_{(2,3)}(u, \tau_1)+Z^5_{(6,1)}(u, \tau_1))\overline{(Z^5_{(3,2)}(v, \tau_2)+Z^5_{(1,6)}(v, \tau_2))},\\
&\\
&|Z_{(1,1)}^9+Z_{(10,1)}^9+Z_{(1,10)}^9|^2+|Z_{(3,3)}^9+Z_{(3,6)}^9+Z_{(6,3)}^9|^2+2|Z_{(4,4)}^9|^2+|Z_{(5,5)}^9+Z_{(5,2)}^9+Z_{(2,5)}^9|^2\\
 &+\left(Z_{(2,2)}^9(u,\tau_1)+Z_{(2,8)}^9(u,\tau_1)+Z_{(8,2)}^9(u,\tau_1)\right)\overline{Z_{(4,4)}^9(v,\tau_2)}+Z_{(4,4)}^9(u,\tau_1)\overline{\left(Z_{(2,2)}^9(v,\tau_2)+Z_{(2,8)}^9(v,\tau_2)+Z_{(8,2)}^9(v,\tau_2)\right)}\\
 &+\left(Z_{(1,4)}^9(u,\tau_1)+Z_{(7,1)}^9(u,\tau_1)+Z_{(4,7)}^9(u,\tau_1)\right)\overline{\left(Z_{(4,1)}^9(v,\tau_2)+Z_{(1,7)}^9(v,\tau_2)+Z_{(7,4)}^9(v,\tau_2)\right)}\\
 &+\left(Z_{(4,1)}^9(u,\tau_1)+Z_{(1,7)}^9(u,\tau_1)+Z_{(7,4)}^9(u,\tau_1)\right)\overline{\left(Z_{(1,4)}^9(v,\tau_2)+Z_{(7,1)}^9(v,\tau_2)+Z_{(4,7)}^9(v,\tau_2)\right)},
 \end{align*}}
respectively. Since the modular invariant of type $E_5^C$ contains
$Z^5_{(1,3)}(u, \tau_1)\overline{Z^5_{(3,1)}(v, \tau_2)},$ but does not contain the term $|Z^5_{(1,3)}|^2$,
 it follows immediately from Lemma \ref{cre} that the modular invariant of type $E_5^C$ cannot be realized by extension of $L_{sl_3}(5,0)$.
Similarly, we can prove that the modular invariants of types $E_9^{(2)}$, $(E_9^{(2)})^C$ cannot be realized by extensions of $L_{sl_3}(9,0)$.\qed
\section{classification of extensions of $L_{sl_3}(k,0)$}
\def\theequation{4.\arabic{equation}}
\setcounter{equation}{0}
\subsection{Extensions of $L_{sl_3}(k,0)$ corresponding to the modular invariant of type $D_k$}
In this subsection we shall show that if $k$ is a positive integer such that $k\equiv 0~(mod~ 3)$, then there exists a unique extension of $L_{sl_3}(k,0)$ corresponding to the modular invariant of type $D_k$. By taking $V=L_{sl_3}(k,0)$, $L=\mathbb{Z}h^{(1)}\oplus \mathbb{Z}h^{(2)}$,  where $h^{(1)}$, $h^{(2)}$ are the fundamental coweights of $sl_3$, we have the following result which was proved in Example 5.11 of \cite{DLM4} and Proposition 3.8, Corollary 3.21 of \cite{L2}.
\begin{theorem}\label{simplecurrent}
Let $k$ be a positive integer such that $k\equiv \text{0} \ (\text{mod}\ 3)$. Then there exists a vertex operator algebra structure on $L_{sl_3}(k, 0)\oplus L_{sl_3}(k, k\Lambda_1)\oplus L_{sl_3}(k, k\Lambda_2)$ such that  $L_{sl_3}(k, 0)\oplus L_{sl_3}(k, k\Lambda_1)\oplus L_{sl_3}(k, k\Lambda_2)$ is an extension of $L_{sl_3}(k, 0)$.
\end{theorem}

 Note that $L_{sl_3}(k, k\Lambda_1)$, $L_{sl_3}(k, k\Lambda_2)$ are simple current $L_{sl_3}(k, 0)$-modules (cf. \cite{DLM4}). By Proposition 5.3 of \cite{DM1}, we have:
  \begin{theorem}\label{unique1}
   There exists a unique vertex operator algebra structure on $L_{sl_3}(k, 0)\oplus L_{sl_3}(k, k\Lambda_1)\oplus L_{sl_3}(k, k\Lambda_2)$ such that  $L_{sl_3}(k, 0)\oplus L_{sl_3}(k, k\Lambda_1)\oplus L_{sl_3}(k, k\Lambda_2)$ is an extension of $L_{sl_3}(k, 0)$.
  \end{theorem}

Moreover, by Proposition 5.2 of \cite{ABD} and  Theorem 3.4 of \cite{HKL}, we know that $L_{sl_3}(k, 0)\oplus L_{sl_3}(k, k\Lambda_1)\oplus L_{sl_3}(k, k\Lambda_2)$ is strongly regular.
Hence, by Theorem \ref{ccc}, $L_{sl_3}(k, 0)\oplus L_{sl_3}(k, k\Lambda_1)\oplus L_{sl_3}(k, k\Lambda_2)$ gives rise to a modular invariant of $L_{sl_3}(k,0)$.
\begin{lemma}\label{modular1}
Let $k$ be a positive integer such that $k\equiv 0 \ (\rm mod\ 3)$. Then the modular invariant corresponding to $L_{sl_3}(k, 0)\oplus L_{sl_3}(k, k\Lambda_1)\oplus L_{sl_3}(k, k\Lambda_2)$ is equal to the the modular invariant of type $D_k$:
$$\frac{1}{3}\sum\limits_{\substack {(m, n)\in P^k \\m\equiv n~(\text{mod}\ 3)}}| Z_{(m, n)}^k + Z_{\sigma(m, n)}^k +Z_{\sigma^2(m, n)}^k|^2.$$
\end{lemma}
\pf By Lemma \ref{5aa},  the modular invariant corresponding to $L_{sl_3}(k, 0)\oplus L_{sl_3}(k, k\Lambda_1)\oplus L_{sl_3}(k, k\Lambda_2)$ must be equal to the modular invariant of type $A_k^C$, $D_k$ or $D_k^C$.
 Note that the modular invariant of type $A_k^C$ is equal to
 \begin{align*}
 \sum\limits_{(m, n)\in P^k} Z_{(m,n)}^k(u,\tau_1)\overline{Z_{(n,m)}^k(v,\tau_2)}.
 \end{align*}
 It follows immediately from Lemma \ref{cre} that the modular invariant of type $A_k^C$ cannot be realized by extension of $L_{sl_3}(k, 0)$.

We next prove that if $k\geq 9$, then the modular invariant of type $D_k^C$ cannot be realized by extension of $L_{sl_3}(k, 0)$. Note that the modular invariant of type $\ D_k^C$ is equal to
\begin{align*}
\frac{1}{3}\sum\limits_{\substack {(m, n)\in P^k \\m\equiv n~(\rm {mod}\ 3)}}&\left( Z_{(m, n)}^k(u,\tau_1)+ Z_{\sigma(m, n)}^k(u,\tau_1) +Z_{\sigma^2(m, n)}^k(u,\tau_1)\right)\\
&\cdot\overline{\left( Z_{h(m, n)}^k(v, \tau_2) + Z_{h\sigma(m, n)}^k(v, \tau_2) +Z_{h\sigma^2(m, n)}^k(v, \tau_2)\right)}.
\end{align*}
By the definition of the maps $h$ and $\sigma$, we have $h(m,n)=(n,m)$, $h\sigma(m,n)=(m, k+3-m-n)$ and $h\sigma^2(m,n)=(k+3-m-n, n)$. If $k\geq 9$, then we have $(m,n)=(1,4)\in P^k$ and  satisfies the property $m\equiv n~(mod ~3)$, $k+3-m-n>4$. Hence, if $k\geq 9$, we then know that the modular invariant of type $D_k^C$ contains the term $Z^k_{(1,4)}(u,\tau_1)\overline{Z^k_{(4,1)}(v,\tau_2)}$, but does not contain the term $|Z^k_{(1,4)}|^2$. It follows from Lemma \ref{cre} that the modular invariant of type $D_k^C$ cannot be realized by extension of $L_{sl_3}(k, 0)$. This forces that the desired modular invariant is equal to the modular invariant of type $D_k$. The proof is complete.
\qed
\subsection{Extension of $L_{sl_3}(k,0)$ corresponding to the modular invariant of type $E_5$}
In this subsection we shall show that there exists a unique extension of $L_{sl_3}(k,0)$ which realizes the modular invariant of type $E_5$. First, it is well-known that there exists a conformal embedding $L_{sl_3}(5,0)\subseteq L_{sl_6}(1,0)$ \cite{SW}.   Moreover, we have the following:
\begin{theorem}\label{2aa}
The affine vertex operator algebra $L_{sl_6}(1,0)$ is an extension of $L_{sl_3}(5,0)$. Furthermore,  $L_{sl_6}(1,0)$ viewed as an $L_{sl_3}(5,0)$-module has the following decomposition
$$L_{sl_6}(1,0)\cong L_{sl_3}(5,0)\oplus L_{sl_3}(5, 2\Lambda_1+2\Lambda_2).$$
\end{theorem}
\pf Note that $L_{sl_3}(5,0)$ and $L_{sl_6}(1,0)$ have the same cental charge \cite{SW}. It follows from Propositions 11.12, 12.12  of \cite{K} that $L_{sl_3}(5,0)$ and $L_{sl_6}(1,0)$ have the same Virasoro vector. Hence, $L_{sl_6}(1,0)$ is an extension of $L_{sl_3}(5,0)$.

We next determine the decomposition of $L_{sl_6}(1,0)$. Since $L_{sl_3}(5,0)$ is rational, $L_{sl_6}(1,0)$ viewed as an $L_{sl_3}(5,0)$-module is a direct sum of finitely many irreducible $L_{sl_3}(5,0)$-modules. Recall that the conformal weight of $L_{sl_3}(5,\lambda)$  is equal to $\frac{\la\lambda,\lambda+2\Lambda_1+2\Lambda_2\ra}{16}$ (see \cite{FZ}). By a direct calculation,  if $M$ is an irreducible $L_{sl_3}(5,0)$-module with integral conformal weights, then  we have $M$ is isomorphic to  $L_{sl_3}(5,0)$ or $L_{sl_3}(5, 2\Lambda_1+2\Lambda_2)$. It follows that $$L_{sl_6}(1,0)\cong L_{sl_3}(5,0)\oplus mL_{sl_3}(5, 2\Lambda_1+2\Lambda_2),$$ for some positive integer $m$. Note that $\dim  L_{sl_3}( 2\Lambda_1+2\Lambda_2)=27$. This implies that
$$L_{sl_6}(1,0)\cong L_{sl_3}(5,0)\oplus L_{sl_3}(5, 2\Lambda_1+2\Lambda_2),$$as desired.
\qed

 Since $L_{sl_6}(1,0)$ is strongly regular, it gives rise to a modular invariant of $L_{sl_3}(5,0)$.
 \begin{lemma}\label{modular3}
 The modular invariant corresponding to $L_{sl_6}(1,0)$ is equal to the modular invariant of type $E_5$:
 \begin{align*}
 &| Z_{(1,1)}^5+ Z_{(3, 3)}^5|^2 +|Z_{(1, 3)}^5+ Z_{(4, 3)}^5|^2 +|Z_{(3, 1)}^5+ Z_{(3 ,4)}^5|^2\\
 &+|Z_{(3, 2)}^5+Z_{(1, 6)}^5|^2+ |Z_{(4, 1)}^5+ Z_{(1, 4)}^5|^2+|Z_{(2, 3)}^5+Z_{(6, 1)}^5|^2.
 \end{align*}
 \end{lemma}
 \pf By Lemmas \ref{5aa}, \ref{modular1}, we know that the modular invariant corresponding to $L_{sl_6}(1,0)$ must be equal to the modular invariant of type $E_5$, $D_k$ $(k\not\equiv  0\ ({\rm mod~ 3)}, ~k\geqq 4)$ or $D_k^C$ $(k\not\equiv  0\ ({\rm mod~ 3)}, ~k\geqq 4)$. Note that the modular invariant of type $D_k^C$ $(k\not\equiv  0\ ({\rm mod~ 3)}, ~k\geqq 4)$ is equal to
 \begin{align*}
 \sum\limits_{(m, n)\in P^k}Z_{(m,n)}^k(u, \tau_1)\overline{Z_{h\sigma^{k(m-n)}(m,n)}^k(v,\tau_2)}.
 \end{align*}
 If $k\geq 4$, we then have $(m,n)=(1,2)\in P^k$ and satisfies the property $k+3-m-n\geq 4$. It follows immediately from Lemma \ref{cre} that the modular invariant of type $D_k^C$ $(k\not\equiv  0\ ({\rm mod~ 3)}, ~k\geqq 4)$ cannot be realized by extension of $L_{sl_3}(k,0)$. Similarly, we can prove that the modular invariant of type $D_k$ $(k\not\equiv  0\ ({\rm mod~ 3)}, ~k\geqq 4)$ cannot be realized by extension of $L_{sl_3}(k,0)$. This forces that the desired modular invariant is equal to the modular invariant of type $E_5$. The proof is complete.
 \qed

We next show that there exists a unique extension of $L_{sl_3}(5,0)$ corresponding to the modular invariant of type $E_5$.
\begin{theorem}\label{6aa}
Let $U$ be an extension  of $L_{sl_3}(5,0)$ such that $U$ is strongly regular and that $U$ viewed as a module of $L_{sl_3}(5,0)$ has the following decomposition
$$U \cong L_{sl_3}(5,0)\oplus L_{sl_3}(5, 2\Lambda_1+2\Lambda_2).$$
Then $U$ is isomorphic to $L_{sl_6}(1,0)$.
\end{theorem}
\pf Let $U$ be an extension of $L_{sl_3}(5,0)$ satisfying the assumption. For any $u\in U$, denote the vertex operator by $Y(u, z)=\sum_{n\in \Z}u_nz^{-n-1}$. Then we know that $U$ is generated by $U_1$ and that $\dim U_1=35$. Moreover, $U_1$ equipped with the bracket $[,]$, which is defined by $[u, v]=u_0v$ for any $u, v\in U_1$, is a reductive Lie algebra  \cite{DM1}. It follows from Corollary 5.9 in \cite{DM2} that $U\cong L_{\g_1}(k_1,0)\otimes \cdots \ot L_{\g_s}(k_s,0)$ for some simple Lie algebras $\g_1, \cdots, \g_s$ and positive integers $k_1, \cdots, k_s$. We next show that $s=1$, that is, $U_1$ is a simple Lie algebra. Otherwise, note that $\g_s$ is a  proper ideal of $U_1$ such that $[\g_s, \g_1\oplus  \cdots\oplus  \g_{s-1}]=0$. Note also that, by assumption, $\g_1\oplus  \cdots\oplus  \g_s$ is a module of $sl_3$, it follows that $\g_s$ viewed as an $sl_3$-module is isomorphic to a proper submodule of
$sl_3\oplus L_{sl_3}(2\Lambda_1+2\Lambda_2).$
However, note that for any proper $sl_3$-submodule $N$ of
$sl_3\oplus L_{sl_3}(2\Lambda_1+2\Lambda_2),$
we have $[sl_3, N]\neq 0$. This is a contradiction. Hence, we have $s=1$. Since $\dim U_1=35$, it follows that $ U_1\cong  sl_6$. This implies that $U\cong L_{sl_6}(1,0)$. The proof is complete.\qed
\subsection{Extension of $L_{sl_3}(k,0)$ corresponding to the modular invariant of type $E_9^{(1)}$} In this subsection we shall show that there exists a unique extension  of $L_{sl_3}(k,0)$ corresponding to the modular invariant of type $E_9^{(1)}$. Consider the affine vertex operator algebra $L_{E_6}(1,0)$. It was shown in \cite{KS} that $\hat{E_6}$ has a conformal subalgebra $\hat{sl_3}$ and that $L_{E_6}(1,0)$ viewed a module of $\hat{sl_3}$ has the following decomposition
\begin{align*}
\begin{split}
L_{E_6}(1,0) \cong& L_{sl_3}(9,0)\oplus L_{sl_3}(9,9\Lambda_2) \oplus L_{sl_3}(9,9\Lambda_1) \oplus L_{sl_3}(9,\Lambda_1+4\Lambda_2)\\
&\oplus L_{sl_3}(9,4\Lambda_1+\Lambda_2)\oplus L_{sl_3}(9,4\Lambda_1+4\Lambda_2).
\end{split}
\end{align*}
Hence, by the similar argument in the proof of Theorem \ref{2aa}, we have:
\begin{theorem}\label{3aa}
The affine vertex operator algebra $L_{E_6}(1,0)$ is an extension of $L_{sl_3}(9,0)$. Moreover, $L_{E_6}(1,0)$ viewed as a module of $L_{sl_3}(9,0)$ has the following decomposition
\begin{align*}
\begin{split}
L_{E_6}(1,0) \cong& L_{sl_3}(9,0)\oplus L_{sl_3}(9,9\Lambda_2) \oplus L_{sl_3}(9,9\Lambda_1) \oplus L_{sl_3}(9,\Lambda_1+4\Lambda_2)\\
&\oplus L_{sl_3}(9,4\Lambda_1+\Lambda_2)\oplus L_{sl_3}(9,4\Lambda_1+4\Lambda_2).
\end{split}
\end{align*}
\end{theorem}
Since $L_{E_6}(1,0)$ is strongly regular, it gives rise to a modular invariant of $L_{sl_3}(9,0)$. By Lemmas \ref{modular2}, \ref{5aa} and Theorem \ref{unique1}, this modular invariant must be equal to the modular invariant of type $E_9^{(1)}$:
 $$|Z_{(1, 1)}^9+ Z_{(1, 10)}^9+Z_{(10, 1)}^9+Z_{(5, 5)}^9+Z_{(5, 2)}^9+Z_{(2, 5)}^9|^2 +2|Z_{(3,3)}^9+Z_{(3,6)}^9+Z_{(6,3)}^9|^2.$$

We next show that there exists a unique extension  of $L_{sl_3}(k,0)$ corresponding to the modular invariant of type $E_9^{(1)}$.
\begin{theorem}\label{8aa}
Let $U$ be an extension of $L_{sl_3}(9,0)$ such that $U$ is strongly regular and that $U$ viewed as a module of $L_{sl_3}(9,0)$ has the following decomposition
\begin{align*}
\begin{split}
U \cong& L_{sl_3}(9,0)\oplus L_{sl_3}(9,9\Lambda_2) \oplus L_{sl_3}(9,9\Lambda_1) \oplus L_{sl_3}(9,\Lambda_1+4\Lambda_2)\\
&\oplus L_{sl_3}(9,4\Lambda_1+\Lambda_2)\oplus L_{sl_3}(9,4\Lambda_1+4\Lambda_2).
\end{split}
\end{align*}
Then $U$ is isomorphic to $L_{E_6}(1,0)$.
\end{theorem}
\pf Let $U$ be an extension of $L_{sl_3}(9,0)$ satisfying the assumption.  For any $u\in U$, denote the vertex operator by $Y(u, z)=\sum_{n\in \Z}u_nz^{-n-1}$. Then we know that  $U_1$ equipped with the bracket $[,]$, which is defined by $[u, v]=u_0v$ for any $u, v\in U_1$, is a reductive Lie algebra  \cite{DM1}. In particular, $U_1\cong \h \oplus \g_1\oplus \cdots \oplus \g_s$ for some simple Lie algebras $\g_1,\cdots, \g_s$ and an abelian Lie algebra $\h$. We first prove that $\h=0$. Otherwise, $\h$ is an $sl_3$-submodule of $U_1$. This implies that $\h$ is isomorphic to a submodule of $L_{sl_3}(\Lambda_1+4\Lambda_2)\oplus L_{sl_3}(4\Lambda_1+\Lambda_2)$. However, for any $sl_3$-submodule $N$ of $L_{sl_3}(\Lambda_1+4\Lambda_2)\oplus L_{sl_3}(4\Lambda_1+\Lambda_2)$, we have $[sl_3, N]\neq 0$. This is a contradiction. By the similar argument in the proof of Theorem \ref{6aa}, we can also prove that $s=1$. In particular, $U_1$ is a simple Lie algebra with dimension $78$.  Hence $U_1$ is isomorphic to $B_6, C_6, \text{or}\ E_6$.

We now let $\tilde U$ be the vertex operator subalgebra of $U$ generated by $U_1$. By Corollary 4.3 of \cite{DM2}, we know that $\tilde U\cong L_{\g_1}(k_1,0)$ for some positive integer $k_1$.  We next show that the Virasoro vector of $\tilde U$ is equal to that of $L_{sl_3}(9,0)$. Otherwise, denote the Virasoro vectors of $\tilde U$, $L_{sl_3}(9,0)$ by $\w'$, $\w$, respectively, then we have $0\neq \w''=\w'-\w$. Moreover, it is known that $[\w''_m, \w_n ]$ for any $m,n\in\Z$ (see Proposition 12.10 of \cite{K}). We then have $\w'_1\w=\w''_1\w+\w_1\w=2\w$ and $\w'_2\w=\w''_2\w+\w_2\w=0$. Therefore, we know that $C_{\tilde U}(L_{sl_3}(9,0))=\{u\in \tilde U|\w_0u=0\}$ is a vertex operator subalgebra of $\tilde U$ and the Virasoro vector of $C_{\tilde U}(L_{sl_3}(9,0))$ is $\w''$ (see \cite{FZ}, \cite{LL}). In particular, $C_{\tilde U}(L_{sl_3}(9,0))$ is infinite dimensional. As a result, we know that $\tilde U$ viewed as an $L_{sl_3}(9,0)$-module is a direct sum of infinitely many irreducible $L_{sl_3}(9,0)$-module (see (12.12.3) of \cite{K}). However, since $L_{sl_3}(9,0)$ is strongly regular, $\tilde U$ viewed as an $L_{sl_3}(9,0)$-submodule of $U$ is a direct sum of finitely many irreducible $L_{sl_3}(9,0)$-modules. This is a contradiction.
Thus, the central charge of $\tilde U$ should be equal to 6. By a direct calculation, we know that the central charges of $L_{B_6}(k_1,0)$
 and $L_{C_6}(k_1,0)$ cannot be equal to $6$.   Therefore, $\tilde U$ must be isomorphic to $L_{E_6}(k_1,0)$. By a direct calculation, we get $k_1=1.$ Thus, $\tilde U\cong L_{E_6}(1,0)$. This implies that $U=\tilde U\cong L_{E_6}(1,0)$. The proof is complete. \qed
\subsection{Extension  of $L_{sl_3}(k,0)$ corresponding to the modular invariant of type $E_{21}$} In this subsection we shall show that there exists a unique  extension  of $L_{sl_3}(k,0)$ corresponding to the modular invariant of type $E_{21}$. Consider the affine vertex operator algebra $L_{E_7}(1,0)$. It was shown in \cite{KS} that $\hat{E_7}$ has a conformal subalgebra $\hat{sl_3}$ and that $L_{E_7}(1,0)$ viewed a module of $\hat{sl_3}$ is isomorphic to
{\small\begin{align}\label{133}
\begin{split}
 & L_{sl_3}(21,0)\oplus L_{sl_3}(21,21\Lambda_1) \oplus L_{sl_3}(21,21\Lambda_2)\oplus L_{sl_3}(21,\Lambda_1+10\Lambda_2)\oplus L_{sl_3}(21,10\Lambda_1+\Lambda_2)\\
 &\oplus L_{sl_3}(21,10\Lambda_1+10\Lambda_2)\oplus L_{sl_3}(21,4\Lambda_1+4\Lambda_2) \oplus L_{sl_3}(21,6\Lambda_1+6\Lambda_2)\oplus L_{sl_3}(21,13\Lambda_1+4\Lambda_2)\\
 &\oplus L_{sl_3}(21,9\Lambda_1+6\Lambda_2) \oplus L_{sl_3}(21,4\Lambda_1+13\Lambda_2) \oplus L_{sl_3}(21,6\Lambda_1+9\Lambda_2).
\end{split}
\end{align}}
Hence, by the similar argument in the proof of Theorem \ref{2aa}, we know that $L_{E_7}(1,0)$ is an extension  of $L_{sl_3}(21, 0)$ and strongly regular.
Therefore, $L_{E_7}(1,0)$ gives rise to a modular invariant of $L_{sl_3}(21,0)$. By Lemma \ref{modular2} and Theorem \ref{unique1}, we know that this modular invariant must be equal to the modular invariant of type $E_{21}$. Moreover, we have the following:
\begin{theorem}\label{4aa}
Let $U$ be an extension  of $L_{sl_3}(21,0)$ such that $U$ is strongly regular and $U$ viewed as a module of $L_{sl_3}(21,0)$ has the decomposition (\ref{133}). Then $U$ is isomorphic to $L_{E_7}(1,0)$.
\end{theorem}
\pf Let $U$ be an extension of $L_{sl_3}(21,0)$ satisfying the assumption. By the similar argument in the proof of Theorem \ref{8aa}, we know that $U_1$ is a simple algebra with dimension 133. It follows that $U_1$ is isomorphic to $E_7$. Let $\tilde U$ be the vertex operator subalgebra of $U$ generated by $U_1$. Then we know that $\tilde U\cong L_{E_7}(k,0)$ for some positive integer $k$ (see Corollary 4.3 of \cite{DM2}).  By the similar argument in the proof of Theorem \ref{8aa}, we can prove that $\tilde U$ and $L_{E_7}(k,0)$ have the same Virasoro vector. In particular, the central charge of $\tilde U$ is equal to $7$. By a direct calculation, we have $k=1$, that is, $\tilde U\cong L_{E_7}(1,0)$.  Hence,  $U$ is isomorphic to $L_{E_7}(1,0)$. The proof is complete.\qed
\subsection{Classification of extensions of $L_{sl_3}(k,0)$} We are now ready to prove the main result in this paper.
\begin{theorem}
Let $U$ be an extension of $L_{sl_3}(k,0)$ such that $U$ is strongly regular. Then $U$ is isomorphic to  $L_{sl_3}(k,0)$,  $L_{sl_6}(1,0)$,  $L_{E_6}(1,0)$,  $L_{E_7}(1,0)$ or the vertex operator algebra in Theorem \ref{simplecurrent}.
\end{theorem}
\pf Let $U$ be a vertex operator algebra satisfying the assumption. By Theorem \ref{ccc}, $U$ gives rise to a modular invariant of $L_{sl_3}(k,0)$. By Lemma \ref{5aa} and the proof of Lemmas \ref{modular1}, \ref{modular3}, we know that the modular invariants of types $A_k^C$,$\ D_k^C\ ~( k\neq 3,6)$, $E_5^C$, $(E_9^{(2)})^C$ and  $D_k ~(k\not\equiv  0~ ( mod~ 3), ~k\geqq 4$), $E_9^{(2)}$ cannot be realized by extensions of $L_{sl_3}(k,0)$. This implies that the modular invariant corresponding to $U$ is equal to the modular invariant of type $A_k$, $D_k$ $(k\equiv 0 ~(mod~ 3))$, $E_5$, $E_9^{(1)}$ or $E_{21}$. Hence, by Theorems \ref{unique1}, \ref{6aa}, \ref{8aa} {\rm and} \ref{4aa}, $U$ is isomorphic to $L_{sl_3}(k,0)$,  $L_{sl_6}(1,0)$,  $L_{E_6}(1,0)$,  $L_{E_7}(1,0)$ or the vertex operator algebra in Theorem \ref{simplecurrent}.
This completes the proof. \qed
\begin{remark}
 It is shown in Lemmas \ref{5aa}, \ref{modular1}, \ref{modular3}, that there exist modular invariants of  $L_{sl_3}(k,0)$, $k\in \Z_+$, cannot be realized by extensions of vertex operator algebras. However, it was shown in \cite{EP} that all the modular invariants of   $L_{sl_3}(k,0)$, $k\in \Z_+$, can be realized by subfactors. Since subfactors are also closely related to conformal field theory, it is interesting to  realize the remain modular invariants in the framework of vertex operator algebras.
\end{remark}
\subsection*{Acknowledgments}
The authors wish to thank the referees for their valuable suggestions.

\end{document}